\documentclass[12pt]{article}
\usepackage{latexsym,amsfonts,amssymb}
\usepackage[dvips]{color}
\usepackage{colordvi,multicol}

\def \cal{\mathcal}

\newtheorem{thm}{Theorem}[section]

\newtheorem{lem}[thm]{Lemma}

\newtheorem{defi}[thm]{Definition}
\newtheorem{rem}[thm]{Remark}

\date{}
\begin{document}

\title{\bf Products of Conditional Expectation Operators: Convergence and Divergence}
 \author{Guolie Lan$^a$, Ze-Chun Hu$^b$  and  Wei Sun$^c$\\ \\ \\
 {\small $^a$   School of Economics and Statistics, Guangzhou University, China}\\
 {\small langl@gzhu.edu.cn}\\ \\
{\small $^b$ College of Mathematics, Sichuan  University,  China}\\
 {\small zchu@scu.edu.cn}\\ \\
 {\small $^c$ Department of Mathematics and Statistics, Concordia
University, Canada}\\
{\small wei.sun@concordia.ca}}

\maketitle



\vskip 0.5cm \noindent{\bf Abstract}\quad In this paper, we
investigate the convergence of  products of
conditional expectation operators. We show that if
$(\Omega,\cal{F},P)$ is a probability space that is not  purely atomic, then
divergent sequences of products of conditional expectation operators involving 3 or 4
sub-$\sigma$-fields of $\cal{F}$ can be constructed for a large
class of random variables in $L^2(\Omega,\cal{F},P)$. This settles
in the negative a long-open conjecture. On the other hand, we show
that if $(\Omega,\cal{F},P)$ is a purely atomic  probability space,
then products of conditional expectation operators involving any finite
set of sub-$\sigma$-fields of $\cal{F}$ must converge for all
random variables  in $L^1(\Omega,\cal{F},P)$.

\smallskip


\smallskip



\section{Introduction and main results}\setcounter{equation}{0}

Conditional expectation is one of the most important concepts in
probability theory. It plays a central role in probability and
statistics. Let $(\Omega,\cal{F},P)$ be a probability space and ${\cal G}_1,{\cal G}_2,\dots,{\cal G}_K$ be sub-$\sigma$-fields of $\cal{F}$, where $K\in \mathbb{N}$. For $k=1,\dots,K$, denote by $E_k$ the conditional expectation operator with respect to $\cal{G}_k$, i.e., $E_k X={E}(X|\,\cal{G}_k)$ for $X\in L^1(\Omega,\cal{F},P)$. Suppose that $\cal{F}_1,\cal{F}_2,\dots\in\{{\cal G}_1,{\cal G}_2,\dots,{\cal G}_K\}$.
For $X_0\in L^1(\Omega,\cal{F},P)$, define the sequence
$\{X_n\}$ successively by
\begin{eqnarray}\label{eqn-int11}
X_n={E}(X_{n-1}|\,\cal{F}_n),\quad n\geq1.
\end{eqnarray}
Then $X_n=E_{k_n} \cdots E_{k_1} X_0$ for some sequence $k_1,k_2,\dots\in \{1,2,\dots,K\}$. In this paper, we will
investigate the convergence of $\{X_n\}$.

Note that conditional expectation operators can be
regarded  as contraction operators on Banach spaces.
The study on the convergence of $\{X_n\}$ is not only of intrinsic interest,
but is also important in various applications including numerical solutions of linear equations and partial differential equations \cite{Browder58,Lion88}, linear inequalities \cite{Spin87},
approximation theory \cite{Reich,Fran84} and computer tomography \cite{Smith77}.

In 1961, Burkholder and Chow \cite{BC61} initiated the study of
convergence of
products  of conditional expectations. They focused on the case $K=2$ and showed that $\{X_n\}_{}$
converges almost everywhere and in $L^2$-norm for $X_0\in
L^2(\Omega,\cal{F},P)$. 
Further, it follows
from Stein \cite{St61} and  Rota \cite{Ro62} that if $X_0\in
L^p(\Omega,\cal{F},P)$ for some $p>1$, then $\{X_n\}_{}$ converges
almost everywhere. On the other hand, Burkholder \cite{Bu62} and
Ornstein \cite{Or68} showed that for $X_0\in
L^1(\Omega,\cal{F},P)$ almost everywhere convergence need not hold
necessarily.

However, for $K\geq 3$ the convergence of
$\{X_n\}_{}$ becomes a very challenging problem.  This paper is
devoted to the following long-open conjecture on convergence of
products  of conditional expectations:

\noindent \textbf{(CPCE)}\ \  If $X_0\in L^2(\Omega,\cal{F},P)$
and all the $\cal{F}_n$ come from a finite set of
sub-$\sigma$-fields of ${\cal F}$, then $\{X_n\}_{}$ must converge
in $L^2$-norm.


Conjecture (CPCE) is
closely related to the convergence of products of orthogonal
projections in Hilbert spaces. Before stating the main results of
this paper, let us recall the important results obtained so far
for the convergence of products of orthogonal projections in
Hilbert spaces.

Let $H$ be a Hilbert space and $H_1,H_2,\dots,H_K$ be
closed subspaces of $H$, where $K\in \mathbb{N}$.  Denote by
$P_{H_k}$ the orthogonal projection of $H$ onto $H_k$. Let $x_0\in
H$ and $k_1,k_2,\dots\in \{1,2,\dots,K\}$, we define the sequence
$\{x_n\}$ by
\begin{eqnarray}\label{eqn-int1}
x_n=P_{H_{k_n}}x_{n-1},\quad n\geq1.
\end{eqnarray}
If $K=2$, the convergence of $\{x_n\}$ in $H$ follows from a
classical result of von Neumann \cite[Lemma 22]{Von}.  If $K\ge 3$
and $H$ is finite dimensional, the convergence of $\{x_n\}$  was
proved by Pr\'ager \cite{Pr}. If $H$ is infinite dimensional and
$\{k_n\}$ is periodic, the convergence of $\{x_n\}$ was obtained
by Halperin \cite{Ha62}. Halperin's result was then generalized to
the quasi-periodic case by Sakai \cite{S}. Based on the results on
the convergence of products  of orthogonal projections, Zaharopol
\cite{Za90},  Delyon and Delyon \cite{DD99},  and  Cohen
\cite{Co07} proved that if  $\{\cal{F}_n\}$  is a periodic
sequence with all the $\cal{F}_n$ coming from a finite set of
$\sigma$-fields, then for any $X_0\in L^p(\Omega,\cal{F},P)$ with
$p>1$, the sequence $\{X_n\}$ of the form (\ref{eqn-int11})
converges in $L^p$-norm and almost everywhere.

In 1965, Amemiya and Ando \cite{AA65} considered the more general
convergence  problem when $K\ge 3$ and $\{k_n\}$ is non-periodic.
They showed that for arbitrary  sequence $\{k_n\}$, the sequence
$\{x_n\}$ of the form (\ref{eqn-int1}) converges weakly in  $H$,
and they posed the question if $\{x_n\}$ converges also in the
norm of $H$. In 2012, Paszkiewicz \cite{Pa12} constructed an
ingenious example of 5 subspaces of $H$ and a sequence $\{x_n\}$
of the form (\ref{eqn-int1}) which does not converge in $H$.
Kopeck\'a and M\"uller resolved in \cite{KM14} fully the question
of Amemiya and Ando. They refined Paszkiewicz's construction to get
an example of 3 subspaces of $H$ and a sequence $\{x_n\}$ which
does not converge in $H$. In \cite{Kop17}, Kopeck\'a and
Paszkiewicz considerably simplified  the construction of
\cite{KM14} and obtained improved  results on the divergence of
products  of orthogonal projections in Hilbert spaces.

 Note that a projection on  $L^2(\Omega,\cal{F},P)$ can not necessary be represented as  a conditional expectation operator.
 Thus, counterexamples for the convergence of products  of orthogonal projections in Hilbert
 spaces do not necessarily yield divergent sequences of products  of conditional expectation on probability spaces.
 Let ${\cal B}(\mathbb{R})$ be the Borel $\sigma$-field of $\mathbb{R}$ and $dx$ be the Lebesgue measure.
 In 2017, Komisarski \cite{Ko17} showed that there exist $X_0\in L^1(\mathbb{R})\cap L^2(\mathbb{R})$ and $\{\cal{F}_n\}$
 coming from 5 sub-$\sigma$-fields of ${\cal B}(\mathbb{R})$ such that the sequence $\{X_n\}$ of the form (\ref{eqn-int11})
 diverges in $L^2(\mathbb{R})$. Note that   $(\mathbb{R},{\cal B}(\mathbb{R}), dx)$ is only a $\sigma$-finite measure space
 and the conditional expectation considered in \cite{Ko17} is understood in an extended sense.
 Conjecture (CPCE)  still remains open for probability spaces. We would like to point out that Akcoglu and King \cite{AK96}
 constructed an example of  divergent sequences involving infinitely many  sub-$\sigma$-fields on the interval $[-\frac{1}{2},\frac{1}{2})$.

In this paper, we will show that Conjecture (CPCE) falls if
$(\Omega,\cal{F},P)$ is not a purely atomic  probability space; however,
it holds if $(\Omega,\cal{F},P)$ is a purely atomic  probability
space. More precisely, we will prove the following results.

\begin{thm} \label{the-Gauss}
There exists a sequence $k_1,k_2,\dots \in\{1,2,3\}$ with the following property:

\noindent Suppose that $X_0$ is a  Gaussian random variable on $(\Omega,\mathcal{F},P)$ and there exists a
non-atomic $\sigma$-field $\cal{C}\subset \mathcal{F}$ which is independent of $X_0$.
Then there exist three $\sigma$-fields $\mathcal{C}_1, \mathcal{C}_2, \mathcal{C}_3\subset \mathcal{F}$,
such that the sequence $\{X_n\}_{}$ defined by
$$X_n=\mathbf{E}(X_{n-1}|\,\mathcal{C}_{k_n}),\quad n\geq1$$
does not converge in probability.
\end{thm}

Denote by $\mathcal{P}(\mathbb{R})$ the space of all probability
measures  on $\mathbb{R}$. Then $\mathcal{P}(\mathbb{R})$  becomes
a complete metric space if it is  equipped with the
L\'{e}vy-Prokhorov metric.

\begin{thm} \label{the-convolGauss}
There exists a sequence $k_1,k_2,\dots \in\{1,2,3\}$ with the following property:

\noindent (1) Suppose that $(\Omega,\cal{F},P)$ is a {non-atomic}
probability space. Then there exists a dense subset ${\cal D}$ of
$\mathcal{P}(\mathbb{R})$ such that for any $\mu\in {\cal D}$ we
can find a random variable $X\in L^2(\Omega,\cal{F},P)$  with
distribution $\mu$  and  four $\sigma$-fields
 $\mathcal{C}_0, \mathcal{C}_1, \mathcal{C}_2, \mathcal{C}_3\subset \mathcal{F}$,
such that the sequence $\{X_n\}_{}$ defined by
$$X_0=\mathbf{E}(X|\,\mathcal{C}_{0}),\ \  X_n=\mathbf{E}(X_{n-1}|\,\mathcal{C}_{k_n}), \ n\geq1$$
does not converge in probability.

\noindent (2) Let $(\Omega,\cal{F},P)$ be a probability space that is not  purely atomic.
Then there exist a random variable $X_0\in L^2(\Omega,\cal{F},P)$ and three $\sigma$-fields
 $\mathcal{G}_1, \mathcal{G}_2, \mathcal{G}_3\subset \mathcal{F}$,
such that the sequence $\{X_n\}_{}$ defined by
$$X_n=\mathbf{E}(X_{n-1}|\,\mathcal{G}_{k_n}), \ n\geq1$$
does not converge in probability.
\end{thm}

Denote by ${\cal N}$ the collection of all null sets of
$(\Omega,\cal{F},P)$. For a sub-$\sigma$-filed ${\cal G}$ of
${\cal F}$, we define $\overline{{\cal G}}$ to be the
$\sigma$-field generated by ${\cal G}$ and ${\cal N}$.

\begin{thm}\label{thm-4.1}
Suppose that $(\Omega,\cal{F},P)$  is a purely atomic probability
space, ${\cal G}_1,\dots,{\cal G}_K$ are sub-$\sigma$-fields of $\cal{F}$, and $\cal{F}_1,\cal{F}_2,\dots\in\{{\cal G}_1,\dots,{\cal G}_K\}$. 
Let $X_0\in
L^p(\Omega,\cal{F},P)$ with $p\geq1$ and $\{X_n\}_{}$ be defined
by (\ref{eqn-int11}). Then $\{X_n\}_{}$  converges to some
$X_{\infty}\in L^p$ in $L^p$-norm and almost everywhere. If each ${\cal G}_k$ repeats
infinitely in the sequence $\{{\cal F}_n\}$, then
$$
 X_{\infty}={E}\left(X_{0}\left|\bigcap_{k=1}^{K}\overline{{\cal G}_k}\right.\right).
$$
\end{thm}

The rest of this paper is organized as follows.
 In Section 2, we discuss the linear compatibility under conditional expectations,
 which is essential for our construction of divergent sequences of products  of conditional expectation operators.
In Section 3, we consider divergent sequences of products  of conditional expectation operators on probability spaces that are not purely atomic  and prove Theorems \ref{the-Gauss} and \ref{the-convolGauss}.
In Section 4, we investigate the convergence of products  of conditional expectation operators on purely atomic  probability spaces and prove Theorem \ref{thm-4.1}.


\section{Linear compatibility and deep uncorrelatedness}\setcounter{equation}{0}

Let $(\Omega,\mathcal{F},P)$ be a probability space.  We consider
the  linear compatibility defined by conditional linear equations,
which is closely related to linear regression and optimal
estimation (cf. Rao \cite{Rao}).

\begin{defi}\label{defi-1.1} Two integrable random variables $X,Y$ on $(\Omega,\mathcal{F},P)$ are said to be {linearly  compatible under conditional expectations},
or {linearly  compatible} in short, if there exist $a,b,c,d \in
\mathbb{R}$ such that almost surely,
\begin{eqnarray}\label{comp}
  \mathbf{E}(Y|X) = aX+c,\ \  \mathbf{E}(X|Y) = bY+d.
\end{eqnarray}
\end{defi}

Obviously, if $X$ and $Y$ are independent or perfectly collinear, i.e., $Y=aX+c$, then they are linearly  compatible.  For non-trivial examples,
note that if $(X,Y)$ have a 2-dimensional Gaussian distribution then they are linearly  compatible, and if both $X$ and $Y$ follow two-point distributions, then they must be linearly  compatible.

\begin{lem} \label{lemab01}
Let $X$ and $Y$ be two random variables on
$(\Omega,\mathcal{F},P)$ with $0<{\rm Var}(X)<\infty$ and $0<{\rm
Var}(Y)<\infty$. Suppose that (\ref{comp}) holds. Denote by
$\rho_{XY}$ the correlation coefficient of  $X$ and $Y$, and
denote by $\sigma(X)$ and $\sigma(Y)$ the $\sigma$-fields
generated by $X$ and $Y$, respectively. Then,
\begin{description}
  \item[(i)] $0\leq ab\leq1$ and ${ab}=\rho_{XY}^2$;
  \item[(ii)] $a b=1$ implies that $Y= a X+c$ a.s.;
  \item[(iii)] $a b<1$ implies that ${E}\left[X|\overline{\sigma(X)}\cap\overline{\sigma(Y)}\right]={E}(X)$ and ${E}\left[Y|\overline{\sigma(X)}\cap\overline{\sigma(Y)}\right]=E(Y)$ a.s.;
  \item[(iv)] $ab=0$ implies that $a=b=0$.
\end{description}
\end{lem}
{\bf Proof.} We assume without loss of generality that $E(X)=E(Y)=0$. Then, $c=d=0$ in (\ref{comp}).

\noindent (i) It follows from (\ref{comp}) that ${E}(XY| X)=X{E}(Y|X)=a X^2$ and ${E}(XY| Y)=Y{E}(X|Y)=b Y^2$. Taking expectations, we get
\begin{equation}\label{exp-xy}
    {E}(XY)=a {E}(X^2)=b {E}(Y^2).
\end{equation}
Then $\rho_{XY}^2=ab$, which implies that $0\leq ab\leq1$.

\noindent (ii) is proved by  Rao (\cite[Proposition 2.1]{Rao}), where only
the finiteness of expectations is assumed.

\noindent (iii) We define the operators $T_X$ and $ T_Y$ on $L^2(\Omega,\mathcal{F},P)$ by  $T_X Z={E}(Z| X)$ and $T_Y Z={E}(Z| Y)$ for  $Z\in L^2(\Omega,\mathcal{F},P)$.
By Burkholder and Chow \cite[Theorem 3]{BC61}, we have
\begin{equation}\label{Tx-Ty}
    (T_X T_Y)^n X\rightarrow {E}\left[X|\overline{\sigma(X)}\cap\overline{\sigma(Y)}\right], \ \
    a.s..
\end{equation}
On the other hand, since  (\ref{comp}) holds with $a b<1$,
\begin{equation}\label{Tx-Ty=ab}
    (T_X T_Y)^n X=(ab)^n X\rightarrow 0, \ \ a.s..
\end{equation}
By (\ref{Tx-Ty}) and (\ref{Tx-Ty=ab}), we get
$${E}\left[X|\overline{\sigma(X)}\cap\overline{\sigma(Y)}\right]=0\ a.s..
 $$
 Similarly, we can show that $${E}\left[Y|\overline{\sigma(X)}\cap\overline{\sigma(Y)}\right]=0\ a.s..$$

\noindent (iv) is a direct consequence of (\ref{exp-xy}) since ${E}(X^2)$ and ${E}(Y^2)$ are non-zero. \hfill\fbox

Motivated by Lemma \ref{lemab01} (iv), we introduce the definition of deep uncorrelatedness for two random variables.

\begin{defi}\label{def-orth}
Two integrable random variables $X,Y$ on $(\Omega,\mathcal{F},P)$
are said to be {deeply uncorrelated} if
$$
    {E}(X|Y)=E(X)\ \ {\rm and}\ \ {E}(Y|X)=E(Y).
$$
\end{defi}

\begin{rem}
It is clear that if $X$ and $Y$ are integrable and independent
then they are deeply uncorrelated, and if $X$ and $Y$ have finite
variances and are deeply uncorrelated then they are
uncorrelated, i.e., $\rho_{XY}=0$. The following examples show
that  deeply uncorrelated is equivalent to neither independent
nor uncorrelated.

\noindent (i) Let $(X,Y)$ be a pair of random variables with the
uniform distribution on the unit disc $\{(x,y):x^2+y^2\leq1\}$. It
can be checked that $X$ and $Y$ are deeply uncorrelated but not
independent.

\noindent (ii)  Let $A$ and $B$ be two measurable sets satisfying
$P(A)=P(B)>0$ and  $P(A\cap B)<P(A\cup B)<1$. Define $X={1}_{A\cup
B}$ and $Y={1}_{A}-{1}_{B}$. Note that $X Y=Y$. Then ${E}(X Y)={E}(Y)=0$, which
implies that $(X,Y)$ are uncorrelated. However, we have that
$$
\int_{\Omega}{E}(X|Y)\cdot Y^2dP={E}(XY^2)={E}(Y^2)=2[P(A)-P(A\cap B)]>0
$$
and
$$
\int_{\Omega}E(X)\cdot Y^2dP=P(A\cup B){E}(Y^2).
$$
Hence ${E}(X|Y)\neq E(X)$, which implies that $X$ and $Y$ are not deeply uncorrelated.
\end{rem}

We now define linear compatibility and deep uncorrelatedness for  a
family of random variables.
\begin{defi}\label{defi-compat}  \textsl{{(1)}} A family of integrable random variables $X_S=\{X_s\}_{ s\in S}$ on $(\Omega,\mathcal{F},P)$
is said to be {linearly  compatible under conditional
expectations}, or {linearly  compatible} in short, if for any
finite sequence $X_0, X_1, \dots, X_n$ in $ X_S$, there exist
$a_0,a_1,\dots, a_n \in \mathbb{R}$ such that almost surely,
\begin{equation}\label{eqn-X0n}
    {E}\left(X_0| X_1, \dots, X_n\right) =a_0+ \sum_{i=1}^n a_i X_i.
\end{equation}
{\textsl{(2)}} $ X_S$ is called a { deeply uncorrelated family}
if for  any finite sequence $X_0, X_1, \dots, X_n\in X_S$,
$$
    {E}\left(X_0|X_1, \dots, X_n\right) ={E}(X_0).
$$
\end{defi}

\begin{rem}
(i) Let $X\in L^2(\Omega,\mathcal{F},P)$ and ${\cal C}\subset{\cal
F}$. It is well-known that ${E}(X|\mathcal{C})$ provides the
$L^2$-optimal estimation of $X$ given $\mathcal{C}$. Thus
(\ref{eqn-X0n}) implies that the $L^2$-optimal estimation of $X_0$
via $X_1, \dots, X_n$ is consistent with the optimal linear
estimation  via $X_1, \dots, X_n$.

\noindent(ii) An important class of linearly  compatible family is
Gaussian processes. For a Gaussian process  $X_T=\{X_t\}_{ t\in
T}$, every finite collection of  random variables $\{X_0, X_1,
\dots, X_n\}\subset X_T$ has a multivariate normal distribution.
Thus (\ref{eqn-X0n}) holds and therefore $X_T$  is linearly
compatible.

\noindent(iii) Let  $\{X_n\}_{ n\geq0}$ be a { deeply uncorrelated family}
with  ${E} X_n=0$, $\forall n\geq 0$. Define  $Y_n=X_0+X_1+\cdots+X_n$. 
Then  $\{Y_n\}_{ n\geq0}$ is a martingale.
\end{rem}

\begin{lem} \label{lem-compat}
Let $X_S=\{X_s\}_{ s\in S}\subseteq L^2(\Omega,\mathcal{F},P)$ be a {linearly  compatible} family
with  ${E}(X_s)=0$, $\forall s\in S$. Then for any infinite sequence  $\{X_n\}_{ n\geq0}\subset X_S$,
there exists $\{a_n\}_{ n\geq1}\subset \mathbb{R}$ such that almost surely,
\begin{equation}\label{eqn-X0inf}
    {E}(X_0| X_i, i\geq1) =\sum_{i=1}^{\infty} a_i X_i.
\end{equation}
In particular, if $ X_S$ is a {deeply uncorrelated family} then for any $\{X_n\}_{ n\geq0}\subset X_S$,
\begin{equation}\label{eqn-inf=0}
    {E}(X_0| X_i, i\geq1) =0.
\end{equation}
\end{lem}
{\bf Proof.} Since  ${E}(X_n)=0$ for ${ n\geq0}$, it follows from
(\ref{eqn-X0n}) that there exists $\{a_{n,m}:1\le m\le
n,n\in\mathbb{N}\}_{}\subset \mathbb{R}$ such that for each $n$,
$$
    {E}(X_0| X_1, \dots, X_n) =\sum_{i=1}^n a_{n,i} X_i.
$$
By the martingale convergence theorem, we have
$$
    {E}(X_0| X_1, \dots, X_n) \rightarrow{E}(X_0| X_i, i\geq1)\quad \mathrm{in}\  L^2(\Omega,\mathcal{F},P),
$$
which implies that $\sum_{i=1}^n a_{n,i} X_i$ converges in
$L^2(\Omega,\mathcal{F},P)$.

Denote by $Y$ the limit of $\sum_{i=1}^n a_{n,i} X_i$  in
$L^2(\Omega,\mathcal{F},P)$. Then there exists $\{a_n\}_{n\ge
1}\subset \mathbb{R}$ such that
$$
Y=\sum_{i=1}^{\infty} a_i X_i. $$ Thus we obtain
(\ref{eqn-X0inf}). The proof of (\ref{eqn-inf=0}) is similar and
we omit the details. \hfill\fbox

Note that for $X\in L^2(\Omega,\mathcal{F},P)$  and ${\cal
C}\subset{\cal F}$, the conditional expectation
${E}(X|\,\mathcal{C})$ can be regarded as the orthogonal
projection of  $X$ onto the closed subspace
$L^2(\Omega,\mathcal{C},P)$. However, in general, the  orthogonal
projection of $L^2(\Omega,\mathcal{F},P)$ onto a closed linear subspace
can not necessary be represented as a
conditional expectation operator. The following lemma shows that
the linear compatibility ensures the one-to-one correspondence
between conditional expectation operator and  orthogonal
projection.

\begin{lem} \label{lem-subace}
Let $H$ be a closed linear subspace of $L^2(\Omega,\mathcal{F},P)$. Suppose that $H$ is a {linearly  compatible} family
with  ${E}(h)=0$ for any $h\in H$. Then for each closed linear subspace $G\subseteq H$ with countable basis,
there exists a sub-$\sigma$-field $\mathcal{G}$ of $\mathcal{F}$ such that for any  $h\in H$,
\begin{equation}\label{eqn-EhGh}
    {E}(h\,|\, \mathcal{G}) =P_G h.
\end{equation}
\end{lem}
{\bf Proof.} Let  $\{g_i,i\geq 1\}$ be an orthonormal basis of $G$
and define $\mathcal{G}=\sigma(g_i, i\geq1)$. By Lemma
\ref{lem-compat}, for any $h\in H$, there exists
$\{a_i\}_{}\subset \mathbb{R}$ such that
\begin{equation}\label{eqn-sub1}
    {E}(h\,|\, \mathcal{G}) =\sum_{i\geq1} a_i g_i.
\end{equation}
Denote $h_0=h-\sum_{i\geq1} a_i g_i.$ Then for any  $g\in G$, almost surely
$$
{E}(\, g h_0\,|\, \mathcal{G})=g \left[{E}(h\,|\, \mathcal{G})-\sum_{i\geq1} a_i g_i\right]=0.
$$
Hence ${E}(\,h_0 g\,)=0$, which implies that $h_0$ and $g$ are orthogonal in  $L^2(\Omega,\mathcal{F},P)$.
Since $g\in G$ is arbitrary, the right hand side of (\ref{eqn-sub1}) equals the
orthogonal projection of  $h$ onto $G$. Therefore, (\ref{eqn-EhGh}) holds.\hfill\fbox

\section{Divergent sequences on probability spaces that are not purely atomic}\setcounter{equation}{0}

\begin{defi}\label{eqn-atomic}  Let $(\Omega,\cal{F},P)$ be a probability space and $\cal{C}$ be a sub-$\sigma$-field  of $\cal{F}$.

\noindent {{(1)}} A  measurable set $B\in\cal{C}$ is
called $\cal{C}$-{atomic} if  $P(B)>0$ and for any $\cal{C}$-measurable set
$A\subset B$,  it holds that either $P(A)=0$ or $P(A)=P(B)$.

\noindent  {{(2)}} $\cal{C}$ is called {non-atomic} if  it contains no
$\cal{C}$-atomic set, i.e., for each $B\in\cal{C}$ with $P(B)>0$,
there exists a  $\cal{C}$-measurable set $A\subset B$ such that
$0<P(A)<P(B)$.

\noindent  {{(3)}} $\cal{C}$ is called {purely atomic} if it
contains a countable number of $\cal{C}$-{atomic} sets $B_1, B_2,
\dots$ such that
$$
    \sum_{i\geq1} P(B_i)=1, \quad P(B_i\cap B_j)=0, \quad i\neq j.
$$
{\textsl{(4)}} $(\Omega,\cal{F},P)$ is said to be {non-atomic} if
$\cal{F}$ is {non-atomic}. $(\Omega,\cal{F},P)$ is said to be
{purely atomic} if $\cal{F}$ is {purely atomic}.
\end{defi}

\begin{rem}
(i) Note that $\cal{C}$ is  {non-atomic} if it is generated by a random variable whose cumulative distribution function (cdf) is continuous   on $\mathbb{R}$, for example, a continuous random variable.
Moreover, $\cal{C}$ is {non-atomic} if and only if there exists a random variable $X\in \cal{C}$
which has a uniform distribution on $(0,1)$ (cf. \cite[\S 2]{P}).

\noindent(ii)  $\cal{C}$ is {purely atomic} if it is generated by a discrete random variable.
Conversely, if $\cal{C}$ is {purely atomic}, then each $\cal{C}$-measurable random variable
has a  discrete distribution on $\mathbb{R}$.
\end{rem}

We now prove  Theorems \ref{the-Gauss} and \ref{the-convolGauss}, which are stated in \S 1. Our proofs are based on the following remarkable result.
\begin{thm}(Kopeck\'a and Paszkiewicz \cite[Theorem 2.6]{Kop17}))\label{thmKK}
There exists a sequence $k_1,k_2,\dots \in\{1,2,3\}$ with the following property:

\noindent If $H$ is an infinite-dimensional  Hilbert Space and $0\neq w_0\in H$,
then there exist three closed subspaces $G_1, G_2, G_3\subset H$,
such that the sequence of iterates  $\{w_n\}_{ n\geq1}$ defined by
$w_n=P_{G_{k_n}}w_{n-1}$ does not converges in $H$.
\end{thm}

\noindent{\bf Proof of Theorem \ref{the-Gauss}}\ \ Denote by $\gamma^1$ the standard Gaussian measure on $\mathbb{R}$ and denote by $
 \gamma^{\infty}= \gamma^1\times \gamma^1\times\cdots$ the standard Gaussian measure on $\mathbb{R}^{\infty}$. For $u\in (0,1)$, we consider its binary representation:
$$
u=\sum_{i=1}^{\infty}  2^{-i}\cdot u[\,i\,],
$$
where $u[1],u[2],\dots\in \{0,1\}$.
For $k\geq1$, define
$$
h_k(u)= \sum_{i=1}^{\infty}  2^{-i}\cdot u[\,2^{k-1}\cdot(2i-1)\,].
$$
Let $\Psi$ be the map
$$
    (0,1)\rightarrow (0,1)^{\infty}, \quad u\mapsto (h_1(u),h_2(u),\dots).
$$
Denote by $dx$ the Lebesgue measure on $(0,1)$. Then it can be
checked that the image measure of $dx$ under $\Psi$ equals the
infinite product measure $(dx)^{\infty}$ on $(0,1)^{\infty}$. Let
$\Phi$ be the cdf of $\gamma^1$.
Define
$$
g_k(u)=\Phi^{-1}\circ h_k(u).
$$
Let $g$ be the   map
$$
    (0,1)\rightarrow \mathbb{R}^{\infty}, \quad u\mapsto (g_1(u),g_2(u),\dots).
$$
Then the image measure of $dx$ under $g$ equals the standard
Gaussian measure $\gamma^{\infty}$ on $\mathbb{R}^{\infty}$.

Since $\cal{C}$ is a
non-atomic sub-$\sigma$-field which is independent of $X_0$,
there exists a random variable $Y\in\cal{C}$ which has a uniform distribution on $(0,1)$  and is independent of $X_0$. Define
$$
Z_k= g_k(Y),\quad k\geq1.
$$
Set
$$
Z_0=\frac{X_0-E(X_0)}{\sqrt{{\rm Var}(X_0)}}.
$$
 Then
$$(Z_0, Z_1, Z_2, \dots)$$
is a sequence of independent standard Gaussian random variables. Let
$$
 H=\overline{\mathrm{{Span}}(Z_0, Z_1, Z_2, \dots)}
$$
be the closed linear span of $(Z_0, Z_1, Z_2, \dots)$. Then $H$ is an infinite-dimension Gaussian Hilbert space,
i.e., a  Gaussian process which is also a Hilbert subspace of $L^2(\Omega,\mathcal{F},P)$.

We now show that  $H$ is a  {linearly  compatible} family. Take
$u, v_1, \dots, v_n\in H$. Let
$$
 V=\mathrm{{Span}}(v_1, \dots, v_n)
$$
be the linear span of $(v_1, \dots, v_n)$. Then the orthogonal
projection of $u$ onto  $V$ can be written as
$$
    P_V u=\sum_{i=1}^n a_i v_i
$$
for some $a_1,\dots, a_n \in \mathbb{R}$. Define
$$
 u_0=u- P_V u=u-\sum_{i=1}^n a_i v_i.
$$
Then  $u_0$ is orthogonal to $V$ and hence is independent of
$\sigma(V)=\sigma( v_1, \dots, v_n)$, since $(u, v_1, \dots, v_n)$
have a joint Gaussian distribution. Therefore,
$$
    {E}\left(u| v_1, \dots, v_n\right)={E}\left.\left(\,u_0+\sum_{i=1}^n a_i v_i\right| v_1, \dots, v_n\right)= \sum_{i=1}^n a_i v_i,
$$
which implies that  $H$ is {linearly  compatible}.

Applying  Theorem  \ref{thmKK} to the infinite-dimensional  Hilbert space  $H$, we find that there exists a sequence $k_1,k_2,\dots \in\{1,2,3\}$ with the following property:

For  $0\neq x_0\in H$, there exist three closed subspaces $G_1, G_2, G_3\subset H$,
such that the sequence  $\{x_n\}_{}$ defined by
$$
    x_n=P_{G_{k_n}}x_{n-1},\quad n\geq1
$$
does not converges in $L^2$-norm.

\noindent Hence there exist three closed subspaces $H_1, H_2, H_3\subset H$,
such that the sequence  $\{X_n\}_{}$ defined by
$$
    X_n=P_{H_{k_n}}X_{n-1},\quad n\geq1
$$
does not converges in $L^2$-norm.

On the other hand,
by Lemma \ref{lem-subace},  there exist three sub-$\sigma$-fields $\mathcal{C}_1, \mathcal{C}_2, \mathcal{C}_3\subset \mathcal{F}$
such that
$$
    {E}(h\,|\, \mathcal{C}_k) =P_{H_k} h,\quad \forall h\in H.
$$
Therefore,
$$
    X_n={E}(X_{n-1}|\,\mathcal{C}_{k_n}),\quad n\geq 1.
$$

Finally, we show that $\{X_n\}_{}$  does not converge in probability.
Suppose that  $X_n$ converges to some $X_{\infty}$ in probability. Note that
$$
{E} (X_n^4)\leq {E} (X_{n-1}^4)\leq\cdots\leq {E} (X_0^4)<\infty,
$$
which implies that $\{X_n^2\}_{}$ is uniformly integrable.
Therefore $X_n\rightarrow X_{\infty}$  in   $L^2$-norm. We have
arrived at a contradiction.\hfill\fbox

\vskip 0.5cm \noindent{\bf Proof of Theorem
\ref{the-convolGauss} (1).} Since $\cal{F}$ is non-atomic, there
exists a random variable $Z\in\cal{F}$ which has a uniform distribution on
$(0,1)$. Following the first part of the proof of Theorem
\ref{the-Gauss}, we can construct three independent standard
Gaussian  random variables $Y_0, Y_1, Y_2$ on
$(\Omega,\mathcal{F},P)$.

For $\varepsilon>0$, let $\gamma_\varepsilon$ be the Gaussian
measure on $\mathbb{R}$ with mean 0 and variance $\varepsilon^2$.
For $\nu\in\mathcal{P}(\mathbb{R})$, define $\nu\ast
\gamma_\varepsilon$ to be the  convolution of $\nu$ and
$\gamma_\varepsilon$, i.e.,
$$
    \int_{-\infty}^{\infty}f(x)\nu\ast \gamma_\varepsilon(dx)=\int_{-\infty}^{\infty}\int_{-\infty}^{\infty} f(x+y) \frac1{\sqrt{2\pi}\varepsilon} e^{-\frac{x^2}{2\varepsilon^2}}dx\nu(dy),\ \ \forall f\in {\cal B}_b(\mathbb{R}).
$$
Then all the moments of $\nu\ast \gamma_\varepsilon$ are finite and  $\nu\ast \gamma_\varepsilon\rightarrow \nu$ weakly as $\varepsilon\rightarrow0$. Define
$$
    \mathcal{P}_{\gamma}(\mathbb{R})=\{\nu\ast \gamma_\varepsilon:\nu\in\mathcal{P}(\mathbb{R}),\ \varepsilon>0\}.
$$
$\mathcal{P}_{\gamma}(\mathbb{R})$ is  a dense subset of
$\mathcal{P}(\mathbb{R})$ with respect to the L\'{e}vy-Prokhorov
metric.

For $\mu\in\mathcal{P}_{\gamma}(\mathbb{R})$ with  $\mu=\nu\ast \gamma_\varepsilon$.
Define
$$
G(y):=\sup\{x\in \mathbb{R}:\nu((-\infty,x])<y\},\ \ 0<y<1.
$$
Then $G\circ\Phi(Y_1)$ has the  probability distribution $\nu$, where $\Phi$ is the cdf of a standard Gaussian random variable. Define
$$
X=\varepsilon Y_0+G\circ\Phi(Y_1).
$$
Then $X$ has the probability distribution $\mu$. Write
$$
\mathcal{C}_0=\sigma(Y_0),\quad \mathcal{C}=\sigma(Y_2).
$$
Then  $$X_0={E}(X|\,\mathcal{C}_{0})=\varepsilon Y_0+{E} (G\circ\Phi(Y_1))$$
is a  Gaussian  random variable which is independent of $\mathcal{C}$.
Therefore, by Theorem \ref{the-Gauss} we can find three sub-$\sigma$-fields
 $\mathcal{C}_1, \mathcal{C}_2, \mathcal{C}_3\subset \mathcal{F}$,
such that the sequence  $\{X_n\}_{}$ defined by
$$ X_n={E}(X_{n-1}|\,\mathcal{C}_{k_n}),\ \ n\geq1$$
does not converge in probability. \hfill\fbox

\vskip 0.5cm \noindent{\bf Proof of Theorem
\ref{the-convolGauss} (2).}
Suppose that $(\Omega,\cal{F},P)$ is neither  purely atomic nor non-atomic. Then there exist an $N\in \mathbb{N}\cup\{\infty\}$ and a sequence of atomic sets $\{B_i\}_{i=1}^N\subset \cal{F}$ such that $\Omega\setminus\bigcup_{i=1}^NB_i\in {\cal F}$ is a non-atomic set. Denote
$$
D:=\bigcup_{i=1}^NB_i,\ \ C:=\Omega \setminus D.
$$
Define
$$
\cal{F}_C=\{A\in\cal{F}:A\subseteq C\},
$$
and
$$
P_C(A)=\frac{P(A)}{P(C)},\ \ A\in \cal{F}_C.
$$
Then $(C,\cal{F}_C,P_C)$ is a  non-atomic  probability space.

For a sub-$\sigma$-field $\mathcal{G}$ of  $\cal{F}_C$, define
\begin{equation}\label{eqn-3GD}
\mathcal{G}\uplus D=\mathcal{G}\cup \{G\cup D: G\in \mathcal{G}\}.
\end{equation}
Then $\mathcal{G}\uplus D$ is a  sub-$\sigma$-field of $\cal{F}$.
For a random variable $X$ on $(C,\cal{F}_C,P_C)$,
we can extend it to a  random variable $1_C\cdot X$ on  $(\Omega,\cal{F},P)$ by defining
$(1_C\cdot X)(\omega)=0$ for $\omega\in D$.
We claim that
\begin{equation}\label{eqn-3E1CXc}
    E(1_C\cdot X|\,\mathcal{G}\uplus D)=1_{C}\cdot E_C(X|\,\mathcal{G}),
\end{equation}
where $E_C(\cdot|\,\mathcal{G})$ denotes the
conditional expectation on $(C,\cal{F}_C,P_C)$.
In fact, the right hand side of (\ref{eqn-3E1CXc}) is obviously $\mathcal{G}\uplus D$-measurable.
Thus it is sufficient to show that
\begin{equation}\label{eqn-3int}
    \int_B 1_C\cdot X d P=\int_B 1_{C}\cdot E_C(X|\,\mathcal{G}) d P
\end{equation}
for any $B\in\mathcal{G}\uplus D$.

Note that by (\ref{eqn-3GD}) we have that $B=G\cup D$ or $B=G$ for some $G\in \mathcal{G}$.

\noindent Case 1: Suppose that $B=G\cup D$ for $ G\in \mathcal{G}$.
Then the left hand side of (\ref{eqn-3int}) is
\begin{eqnarray*}
  \int_{G\cup D} 1_C\cdot X d P &=& \int_{G}  X\, d P  \\
   &=& \int_{G}  P(C)\cdot X\, d P_C  \\
   &=& \int_{G}  P(C)\cdot E_C(X|\,\mathcal{G})\, d P_C \\
   &=& \int_{G}  E_C(X|\,\mathcal{G})\, d P \\
   &=& \int_{G\cup D} 1_{C}\cdot E_C(X|\,\mathcal{G})\, d P.\\
\end{eqnarray*}
Thus (\ref{eqn-3int}) holds.

\noindent Case 2:  Suppose that $B=G$ for some $ G\in \mathcal{G}$. The proof is similar to that of Case 1 and we omit the details.

Note that  $(C,\cal{F}_C,P_C)$ is a  non-atomic  probability space.
Then there exists on  $(C,\cal{F}_C,P_C)$ a random variable $Z_C$
 which has a uniform distribution on $(0,1)$.
Following the first part of the proof of Theorem
\ref{the-Gauss}, we can construct on $(C,\cal{F}_C,P_C)$ two independent standard
Gaussian  random variables $Y_{0}$ and $Y_C$.
Let
$$
\mathcal{C}=\sigma(Y_C).
$$
Then  $\cal{C}$ is a
non-atomic sub-$\sigma$-field of $\cal{F}_C$ and is independent of the Gaussian  random variable $Y_{0}$.
Thus, by Theorem \ref{the-Gauss}, we can find three sub-$\sigma$-fields
 $\mathcal{C}_1, \mathcal{C}_2, \mathcal{C}_3$  of $\cal{F}_C$
such that the sequence of iterates  $\{Y_{n}\}_{}$ on  $(C,\cal{F}_C,P_C)$ defined by
\begin{equation}\label{eqn-3Yn}
     Y_{n}={E}_C\,(\,Y_{n-1}|\,\mathcal{C}_{k_n}),\ \ n\geq1
\end{equation}
does not converge in probability.
Hence  $\{Y_{n}\}_{}$ must diverge in the $L^2$-norm of  $(C,\cal{F}_C,P_C)$, i.e.,
\begin{equation}\label{eqn-3sup-YnYm}
    \limsup_{n,m\rightarrow\infty}\|Y_n -Y_m \|_{2,\, C}>0,
\end{equation}
where $\|Y_{}\|_{2,\, C}:=[{E}_C(\,Y_{}^2)]^{1/2}$ is  the $L^2$-norm of  $(C,\cal{F}_C,P_C)$ for $Y\in L^2(C,\cal{F}_C,P_C)$.

Now we can construct on $(\Omega,\cal{F},P)$ three sub-$\sigma$-fields
 $\mathcal{G}_1, \mathcal{G}_2, \mathcal{G}_3$ by
\begin{equation}\label{eqn-3calGi}
     \mathcal{G}_i=\mathcal{C}_i\uplus D, \ i=1,2,3,
\end{equation}
and a sequence of random variables $\{X_n\}_{}$ by
\begin{equation}\label{eqn-3XnYn}
     X_n= 1_C\cdot Y_n, \ \ n\geq0.
\end{equation}

Note that by (\ref{eqn-3E1CXc}), (\ref{eqn-3Yn}), (\ref{eqn-3calGi}) and (\ref{eqn-3XnYn}) we  have
$$
\mathbf{E}(X_{n-1}|\,\mathcal{G}_{k_n})=\mathbf{E}(1_C\cdot Y_{n-1}|\,\mathcal{G}_{k_n})
=1_{C}\cdot {E}_C\,(\,Y_{n-1}|\,\mathcal{C}_{k_n})=1_{C}\cdot Y_{n}=X_{n},
$$
which implies that
$$
     X_n=\mathbf{E}(X_{n-1}|\,\mathcal{G}_{k_n}), \ n\geq1.
$$

To show that  $\{X_n\}_{}$ does not converge in probability, note that
\begin{eqnarray}\label{lan22}
  \|X_n -X_m \|^2_{2} &=& \int_{\Omega} (1_{C}\cdot Y_{n}-1_{C}\cdot Y_{m})^2\, d P  \nonumber\\
   &=& \int_{C} ( Y_{n}-Y_{n})^2\, d P   \nonumber\\
   &=& \int_{C} ( Y_{n}-Y_{m})^2\cdot P(C)\, d P_C   \nonumber\\
   &=& P(C)\cdot \|Y_n -Y_m \|^2_{2,\, C}.
\end{eqnarray}
Then we obtain by (\ref{eqn-3sup-YnYm}) and (\ref{lan22}) that
\begin{equation}\label{eqn-3sup|Xn-Xm}
    \limsup_{n,m\rightarrow\infty}\|X_n -X_m \|_{2}=\sqrt{ P(C)}\cdot \limsup_{n,m\rightarrow\infty}\|Y_n -Y_m \|_{2,\, C}>0.
\end{equation}
On the other  hand, we can check that
$$
{E} (X_n^4)\leq {E} (X_{n-1}^4)\leq\cdots\leq {E} (X_0^4)=P(C)\cdot {E}_C (Y_0^4)<\infty,
$$
which implies that $\{X_n^2\}_{}$ is uniformly integrable.
Suppose that  $\{X_n\}$ converges in probability.
Then  $\{X_{n}\}_{}$ converges in the $L^2$-norm,
which contradicts with (\ref{eqn-3sup|Xn-Xm}).
Therefore, $\{X_n\}_{}$  does not converge in probability. \hfill\fbox

\section{Convergence on purely atomic  probability spaces}\setcounter{equation}{0}

In this section we will prove Theorem \ref{thm-4.1}, which is
stated in \S 1. First, we give a lemma that holds for any
probability space.

\begin{lem}\label{thm-4.4}
Let $(\Omega,\cal{F},P)$ be a probability space,  $\{\cal{F}_n\}$
be a family of  sub-$\sigma$-fields of ${\cal F}$ and $\{X_n\}_{}$
be defined by (\ref{eqn-int11}). Then the following statements are
equivalent.

\noindent\textsl{ (1)}
  $\{X_n\}$ converges in ${L^p}$-norm for any $X_0\in L^p(\Omega,\cal{F},P)$ with $p\ge 1$.

\noindent\textsl{ (2)}
  $\{X_n\}$ converges in probability for any $X_0\in L^{\infty}(\Omega,\cal{F},P)$.

\end{lem}

\noindent {\bf Proof.}  It is sufficient to show that
$(2)\Rightarrow(1)$.

For $X\in L^{p}(\Omega,\cal{F},P)$, we define the operators $T_n$
recursively by $T_0 X=X$ and
$$
T_n X = {E}(T_{n-1}X|\,\cal{F}_n),\quad n\geq1.
$$
Then  $\{T_n\}$ is a family of linear  contraction operators on
$L^{p}(\Omega,\cal{F},P)$ and $\sup_{n\ge
1}\|T_nY\|_{\infty}\le\|Y\|_{\infty}$ for $Y\in
L^{\infty}(\Omega,\cal{F},P)$. Let $X_0\in
L^{p}(\Omega,\cal{F},P)$. Then $T_n X_0=X_n$ for $n\ge 1$.

Suppose that (2) holds, i.e., the sequence $\{T_n Y\}$ converges
in probability for each $Y\in L^{\infty}(\Omega,\cal{F},P)$. Then
we obtain by the bounded convergence theorem that
\begin{eqnarray*}
&&\limsup_{n,m\rightarrow\infty}\|T_n X_0-T_m X_0\|_p\\
&\leq&
\limsup_{k\rightarrow\infty}\limsup_{n,m\rightarrow\infty}\{\|T_n(X_0
1_{\{|X_0|> k\}})\|_p \\
&&\quad \quad \quad \quad \quad \quad
 +\|T_n(X_01_{\{|X_0|\le k\}})-T_m(X_0 1_{\{|X_0|\le k\}}) \|_p
 +\|T_m(X_01_{\{|X_0|> k\}})\|_p\}\\
&\leq&2\lim_{k\rightarrow\infty}\|X_0 1_{\{|X_0|> k\}}\|_p+
\limsup_{k\rightarrow\infty}\limsup_{n,m\rightarrow\infty}\|T_n(X_01_{\{|X_0|\le
k\}})-T_m(X_0
1_{\{|X_0|\le k\}}) \|_p\\
&=&0.
\end{eqnarray*}
Hence  $\{T_nX_0\}$ is a Cauchy sequence in
$L^{p}(\Omega,\cal{F},P)$ and therefore converges in $L^p$-norm.
\hfill\fbox

\vskip 0.5cm

\noindent {\bf Proof of Theorem \ref{thm-4.1}.} Let $B_1, B_2,
\dots \in\cal{F}$ be atomic sets (cf. Definition \ref{eqn-atomic})
such that
$$
    \sum_{i\geq1} P(B_i)=1, \quad P(B_i\cap B_j)=0, \quad i\neq j.
$$
Let $X_0\in L^{\infty}(\Omega,\cal{F},P)$. Then for each $X_n$,
there exists a sequence $b_{n, \,1}, b_{n,\, 2}, \dots \in
\mathbb{R}$ such that
$$
    X_n=\sum_{i\geq1} b_{n,\, i}1_{B_i}\quad a.s..  
$$

We define the orthogonal projections $P_n$, $n\ge 1$, on
$L^2(\Omega,\mathcal{F},P)$ by
$$P_n Y={E}(Y|{\cal F}_n),\ \ Y\in
L^2(\Omega,\mathcal{F},P).$$ Then $ P_n\cdots P_2P_1X_0=X_n$. By
Amemiya-Ando \cite[Theorem]{AA65}, $X_n$ converges weakly in
$L^2(\Omega,\cal{F},P)$. Thus
$b_{n,i}=\frac{1}{P({B_i})}\int_{B_i}X_ndP$ converges as
$n\rightarrow\infty$  for each $i\ge 1$, which implies that $X_n$
converges almost everywhere. Since $X_0\in
L^{\infty}(\Omega,\cal{F},P)$ is arbitrary, we obtain by Lemma
\ref{thm-4.4} that $\{X_n\}$ converges to some $X_{\infty}\in L^p$
in ${L^p}$-norm for any $X_0\in L^p(\Omega,\cal{F},P)$ with $p\ge
1$. Further, $\{X_n\}$ converges to $X_{\infty}$ almost everywhere
since $(\Omega,\cal{F},P)$ is a purely atomic probability space.

We now show that
$X_{\infty}={E}(X_{0}|\bigcap_{k=1}^{K}\overline{{\cal G}_k})$.
For each $k\in \{1,\cdots,K\}$,  we can find an
infinite subsequence
 $k_1, k_2, \dots\in\mathbb{N}$  such that $\cal{F}_{k_n}={\cal G}_k,\ n\geq1$.
 It follows that
$$
X_{k_n}\in\overline{\cal{F}_{k_n}}=\overline{{\cal G}_k}.
$$
By the  almost sure convergence of $\{X_n\}_{}$, we have that
$$
X_{k_n}\stackrel{a.s.}{\longrightarrow} X_{\infty},
$$
which implies that  $X_{\infty}\in\overline{{\cal G}_k}$ for each
${\cal G}_k$. Hence
$$
X_{\infty}\in\textstyle\bigcap_{k=1}^{K} \overline{{\cal G}_k}.
$$

Let $A\in \bigcap_{k=1}^{K}\overline{{\cal G}_k}=\bigcap_{n=1}^{\infty}\overline{\cal{F}_{n}}$. By
(\ref{eqn-int11}), we get
$$
{E}(X_n 1_{A})={E}(X_{n-1} 1_{A})=\cdots= {E}(X_{0} 1_{A}).
$$
Then,
\begin{eqnarray*}
  |{E}(X_01_{A})-{E}(X_{\infty} 1_{A})|
   &=& |{E}(X_n1_{A})-{E}(X_{\infty} 1_{A})|  \nonumber \\
   &\leq& {E}(| X_n 1_{A}-X_{\infty} 1_{A}| ) \nonumber\\
   &\leq&  {E}(| X_n -X_{\infty} |) \nonumber\\
   &\rightarrow&0\ \ {\rm as}\ n\rightarrow\infty.
\end{eqnarray*}
Thus
$$
{E}(X_0 1_{A})={E}(X_{\infty} 1_{A}).
$$
Since $A\in \bigcap_{k=1}^{K}\overline{{\cal G}_k}$ is
arbitrary, the proof is complete.\hfill\fbox

\bigskip

{ \noindent {\bf\large Acknowledgments} \quad   This work was
supported by the China Scholarship Council (No. 201809945013),
National Natural Science Foundation of China (No. 11771309) and
Natural Sciences and Engineering Research Council of Canada.}


\begin{thebibliography}{1234}

\bibitem{AK96} M. Akcoglu, J. L. King: An example of pointwise non-convergence of iterated conditional expectation operators, Israel J. Math., 94, 179-188 (1996).

\bibitem{AA65} I. Amemiya, T. Ando: Convergence of random pruducts of contractions in Hilberet space, Acta Scientiarum Mathematicarum (Szegeed), 26, 239-244 (1965).



\bibitem{P} P. Berti, L. Pratelli, P. Rigo: Atomic intersection of $\sigma$-fields and some of its consequences, Probab. Theory Relat. Fields, 148, 269-283 (2010).

\bibitem{Browder58} F. E. Browder,  On some approximation methods for solutions of the Dirichlet problem for linear elliptic equations of arbitrary order, J. Math. Mech., 7, 69-80 (1958).

\bibitem{Bu62} D. L. Burkholder: Successive conditional expectations of an integrable function, Ann. Math. Stat., 33, 887-893 (1962).


\bibitem{BC61} D. L. Burkholder, Y. S. Chow:  Iterates of conditional expectation operators, Proc. Amer. Math. Soc.,  12, 490-495 (1961).

\bibitem{Co07} G. Cohen: Iterates of a product of conditional expectation operators, J. Func. Anal., 242, 658-668 (2007).

\bibitem{DD99} B. Delyon, F. Delyon: Generalization of von Neumann's spctral sets and integral representation of operators, Bull. Soc. Math. France, 127, 25-41 (1999).


\bibitem{Fran84} C. Franchetti, W. Light: The alternating algorithm in uniformly convex space, J. London Math. Soc., 29, 545-555 (1984).


\bibitem{Ha62} I. Halperin: The product of projection operators, Acta Sci. Math. (Szeged), 23, 96-99 (1962).



\bibitem{Ko17} A. Komisarski: Compositions of conditional expectations, Amemiya-And\^{o} conjecture and paradoxes of thermodynamics, J. Func. Anal., 273, 1195-1119 (2017).

\bibitem{KM14} A. Kopeck\'{a}, V. M\"{u}ller: A product of three projections, Studia Math., 223, 175-186 (2014).

\bibitem{Kop17} A. Kopeck\'{a}, A. Paszkiewicz: Strange products of projections, Israel J. Math., 219, 271-286 (2017).


\bibitem{Lion88} P. L. Lions: On the Schwarz alternating method. I, Domain decomposition methods for partial differential equations, SIAM, Philadelphia,  1-42 (1988).

\bibitem{Von} J. von Neumann: On rings of operators. Reduction theory, Ann. of Math., 50, 401-485 (1949).

\bibitem{Or68} D. Ornstein: On the pointwise behavior of iterates of a self-adjoint operators, J. Math. Mech., 18, 473-477 (1968).


\bibitem{Pa12} A. Paszkiewicz: The Amemiya-Ando conjecture falls, arXiv: 1203.3354 (2012).

\bibitem{Pr} M. Pr\'ager: \"Uber ein Konvergenzprinzip im Hilbertschen Raum, Czechoslovak
Math. J., 10, 271-282 (1960).

\bibitem{Rao} M. M. Rao, Exploring ramifications of the equation $E(Y|X)=X$, J. Stat. Theo. Prac., 1, 73-88 (2007).


\bibitem{Reich} S. Reich, Nonlinear semigroups, accretive operators and applications, { Nonlinear Phenomena in Mathematical Sciences}, Academic Press, New York, 831-838 (1982).


\bibitem{Ro62} G. C. Rota: An ``Alternierende Verfahren" for general positive operators, Bull. Amer. Math. Soc., 68, 95-102 (1962).

\bibitem{S} M. Sakai: Strong convergence of infinite products of orthogonal projections
in Hilbert space, Appl. Anal., 59, 109-120 (1995).

\bibitem{Smith77} K. T. Smith, D. C. Solmon, S. L. Wagner: Practical and mathematical aspects of reconstructing objects from radiographs, Bull. Amer. Math. Soc., 83, 1227-1270 (1977).

\bibitem{Spin87} J. E. Spingarn: A projection method for least-squares solutions to over-determined systems of linear inequalities, Linear Algebra Appl., 86, 211-236 (1987).

\bibitem{St61} E. M. Stein: On the maximal ergodic theorem, Proc. Natl. Acad. Sci. USA, 47, 1894-1897 (1961).


\bibitem{Za90} R. Zaharopol: On products of conditional expectation operators,  Canad. Math. Bull., 33, 257-260 (1990).




\end{thebibliography}
\end{document}